\documentclass[a4paper,12pt]{article}
\usepackage[latin1]{inputenc}
\usepackage{amssymb}
\title{{\bf  The Isomorphism Relation  Between Tree-Automatic  Structures } }
\author{Olivier Finkel$^1$    ~~~~ and  ~~~~  Stevo Todor{\v{c}}evi{\'c}$^{1, 2}$
\\
\\
  $^1${\it Equipe de Logique Math\'ematique }
 \\CNRS and  Universit\'e Paris 7,  France.
\\finkel@logique.jussieu.fr
\\ stevo@logique.jussieu.fr
\\
\\ $^2${\it  Department of Mathematics}
\\ University of Toronto, Toronto, Canada
M5S 2E4.
}

\date{}
\begin{document}

\newtheorem{The}{Theorem}[section]
\newtheorem{Pro}[The]{Proposition}
\newtheorem{Deff}[The]{Definition}
\newtheorem{Lem}[The]{Lemma}
\newtheorem{Rem}[The]{Remark}
\newtheorem{Exa}[The]{Example}
\newtheorem{Cor}[The]{Corollary}
\newtheorem{Notation}[The]{Notation}

\newcommand{\vp}{\varphi}
\newcommand{\lb}{\linebreak}
\newcommand{\fa}{\forall}
\newcommand{\Ga}{\Gamma}
\newcommand{\Gas}{\Gamma^\star}
\newcommand{\Gao}{\Gamma^\omega}
\newcommand{\Si}{\Sigma}
\newcommand{\Sis}{\Sigma^\star}
\newcommand{\Sio}{\Sigma^\omega}
\newcommand{\ra}{\rightarrow}
\newcommand{\hs}{\hspace{12mm}

\noi}
\newcommand{\lra}{\leftrightarrow}
\newcommand{\la}{language}
\newcommand{\ite}{\item}
\newcommand{\Lp}{L(\varphi)}
\newcommand{\abs}{\{a, b\}^\star}
\newcommand{\abcs}{\{a, b, c \}^\star}
\newcommand{\ol}{ $\omega$-language}
\newcommand{\orl}{ $\omega$-regular language}
\newcommand{\om}{\omega}
\newcommand{\nl}{\newline}
\newcommand{\noi}{\noindent}
\newcommand{\tla}{\twoheadleftarrow}
\newcommand{\de}{deterministic }
\newcommand{\proo}{\noi {\bf Proof.} }
\newcommand {\ep}{\hfill $\square$}

\maketitle

\begin{abstract}
\noi   An $\om$-tree-automatic structure is a relational structure whose domain and relations are accepted by Muller or Rabin tree automata.
We investigate in this paper the isomorphism problem for  $\om$-tree-automatic structures. We prove first that the  isomorphism relation for
$\om$-tree-automatic boolean algebras (respectively, partial orders, rings, commutative rings,
non commutative rings, non commutative groups, nilpotent  groups of class $n\geq 2$) is not determined by the  axiomatic system {\bf ZFC}.
Then we prove that
the  isomorphism problem for  $\om$-tree-automatic boolean algebras (respectively, partial orders, rings, commutative rings,
non commutative rings, non commutative groups, nilpotent  groups of class $n\geq 2$)
is neither a $\Si_2^1$-set nor a  $\Pi_2^1$-set.
\end{abstract}

\noi {\small {\bf Keywords:}    $\om$-tree-automatic structures;  boolean algebras;
partial orders; rings; groups; isomorphism relation; models of set theory; independence results.}

\section{Introduction}

An automatic structure is a  relational structure whose domain and relations are recognizable by finite automata reading finite words.  Automatic structures
 have very nice decidability and definability properties and have been much studied in the last few years,
see \cite{BlumensathGraedel04,KNRS,NiesBSL,RubinPhd,RubinBSL}.
They form a subclass of the class of (countable) recursive structures where ``recursive" is replaced by
``recognizable by finite automata".
 Blumensath considered in \cite{Blumensath99} more powerful kinds of automata. If we replace automata by tree automata (respectively, B\"uchi automata
reading infinite words, Muller or Rabin tree automata reading infinite labelled trees) then we get the notion of  tree-automatic (respectively,
$\om$-automatic, $\om$-tree-automatic) structures.
In particular, an $\om$-automatic or $\om$-tree-automatic structure may have uncountable cardinality.
All these kinds of automatic structures have the two following fundamental properties.
$(1)$ The class
of automatic (respectively,   tree-automatic, $\om$-automatic, $\om$-tree-automatic) structures is closed under first-order interpretations.
$(2)$ The first-order theory of an automatic (respectively,   tree-automatic, $\om$-automatic, $\om$-tree-automatic) structure is decidable. 

\hs A natural problem is to classify firstly automatic  structures (presentable by finite automata)
using some invariants. For instance Delhomm\'e proved that  the automatic ordinals are the ordinals smaller
than $\om^\om$,    see \cite{Delhomme,RubinPhd,RubinBSL}.  And   
Khoussainov, Nies, Rubin, and Stephan proved in \cite{KNRS} that 
the automatic infinite boolean algebras are the finite products $B_{fin-cof}^n$ of the boolean algebra
 $B_{fin-cof}$ of finite or cofinite subsets of the set of positive integers $\mathbb{N}$. On the other hand some classes of automatic structures,
like automatic linear orders, or automatic groups, are not completely  determined.
 Another fundamental question which naturally arises in the investigation of   the richness of the class of automatic structures is
the following: ``what is the complexity of the isomorphism problem for the class of automatic structures, or for a subclass of it?"
Khoussainov, Nies, Rubin, and  Stephan proved in  \cite{KNRS} that    the isomorphism problem
for the class of automatic structures, or  even for the class of automatic graphs,
is $\Si_1^1$-complete, i.e. as complicated as the isomorphism
problem for recursive structures.  However for some classes of automatic structures, like the classes of automatic ordinals or of
automatic  boolean algebras, the isomorphism problem is decidable, \cite{KNRS,RubinBSL}.
But for other classes like the classes of automatic linear orders or groups the complexity or even the decidability of the isomorphism problem is still unknown.

\hs There has been less classification work for  $\om$-automatic and $\om$-tree-automatic structures.
In particular, it seems that no complete classification exists for classes of $\om$-automatic or $\om$-tree-automatic structures, like the result classifying
completely the automatic ordinals or the automatic boolean algebras.
Some foundational questions about
$\om$-automatic structures have been recently solved.
Kuske and Lohrey proved in \cite{KuskeLohrey} that the first-order theory,  extended with some cardinality quantifiers,
of an (injectively)  $\om$-automatic structure   is decidable. Next Barany, Kaiser and Rubin extended this result to all $\om$-automatic structures and
proved that an $\om$-automatic structure which is countable  is  automatic, i.e. presentable by automata reading finite words,  \cite{KRB08}.
One of the most important foundational problems in this area  is again the question of the complexity of the isomorphism problem for  $\om$-automatic or
$\om$-tree-automatic  structures.
In a recent paper Hjorth, Khoussainov, Montalb{\'a}n, and Nies proved
that the isomorphism problem for  $\om$-automatic structures is not a   $\Si_2^1$-set, \cite{HjorthKMN08}.
In fact their proof implies also  that this isomorphism problem is not a $\Pi_2^1$-set.
Moreover  this is also the case for the restricted class of  $\om$-automatic (abelian) groups and   for the class of  all $\om$-tree-automatic structures
which is an extension of the class of $\om$-automatic structures.

\hs We investigate in this paper the isomorphism problem for some classes of  $\om$-tree-automatic structures.
We prove first that the  isomorphism relation for $\om$-tree-automatic structures (respectively,  $\om$-tree-automatic boolean algebras)
 is not determined by the  axiomatic system {\rm ZFC}.  Indeed, using known results about quotients of the boolean algebra $\mathcal{P}(\mathbb{N})$ over
 analytic ideals on $\mathbb{N}$, we prove that there exist two
$\om$-tree-automatic boolean algebras $\mathcal{B}_1$ and
$\mathcal{B}_2$ such that: $(1)$ $\mathcal{B}_1$ is isomorphic to
$\mathcal{B}_2$ in {\rm  (ZFC + CH)} and $(2)$    $\mathcal{B}_1$
is not  isomorphic to $\mathcal{B}_2$ in {\rm   (ZFC + OCA)},
where the axioms {\rm  CH},  {\rm   OCA} denote respectively the
Continuum Hypothesis, the Open Coloring Axiom. Then we infer from
this result that  the isomorphism relation for
$\om$-tree-automatic  partial orders (respectively, rings,
commutative rings, non commutative rings, non commutative groups,
nilpotent  groups of class $n\geq 2$) is not determined by the
axiomatic system {\rm ZFC}. This shows the importance of different
axiomatic systems of Set Theory in the area of
$\om$-tree-automatic structures. Then we prove  that the
isomorphism problem for $\om$-tree-automatic boolean algebras
(respectively, partial orders, rings, commutative rings, non
commutative rings, non commutative groups, nilpotent groups of
class $n\geq 2$) is neither a $\Si_2^1$-set nor a $\Pi_2^1$-set.

\hs The paper is organized as follows. In Section 2 we recall definitions and first properties of $\om$-automatic and $\om$-tree-automatic structures. We 
give in Section 3 the two   boolean algebras $\mathcal{B}_1$ and
$\mathcal{B}_2$ and prove that they are $\om$-tree-automatic atomless boolean algebras. In Section 4 we introduce notions of topology and prove some 
topological properties which will be useful in the sequel. We recall some notions of Set Theory in Section 5. We prove our main results in Section 6.

\section{$\om$-tree-automatic structures}

\noi When $\Si$ is a finite alphabet, a {\it non-empty finite word} over $\Si$ is any
sequence $x=a_1\ldots a_k$, where $a_i\in\Sigma$
for $i=1,\ldots ,k$ , and  $k$ is an integer $\geq 1$. The {\it length}
 of $x$ is $k$.
 The {\it empty word} has no letter and is denoted by $\varepsilon$; its length is $0$.
 For $x=a_1\ldots a_k$, we write $x(i)=a_i$. 
 $\Sis$  is the {\it set of finite words} (including the empty word) over $\Sigma$.
 
\hs  The {\it first infinite ordinal} is $\om$.
 An $\om$-{\it word} over $\Si$ is an $\om$ -sequence $a_1 \ldots a_n \ldots$, where for all
integers $ i\geq 1$, ~
$a_i \in\Sigma$.  When $\sigma$ is an $\om$-word over $\Si$, we write
 $\sigma =\sigma(1)\sigma(2)\ldots \sigma(n) \ldots $,  where for all $i$,~ $\sigma(i)\in \Si$.
\nl
 The {\it set of } $\om$-{\it words} over  the alphabet $\Si$ is denoted by $\Si^\om$.
An  $\om$-{\it language} over an alphabet $\Sigma$ is a subset of  $\Si^\om$.  

\hs   We consider in this paper relational structures which are presentable by  automata reading infinite trees or   infinite words.
We  assume that the reader is familiar with the notion of B\"uchi automaton reading
infinite words over a finite alphabet which can be found for instance in \cite{Thomas90,Staiger97}.  Informally speaking an 
$\om$-word $x$ over $\Si$  is accepted by a B\"uchi automaton $\mathcal{A}$ iff there is an infinite run of 
$\mathcal{A}$ on $x$  enterring  infinitely often in  some final state of $\mathcal{A}$. The $\om$-language 
$L(\mathcal{A})\subseteq \Sio$  accepted by the B\"uchi automaton $\mathcal{A}$ is
the set of  $\om$-words $x$ accepted by $\mathcal{A}$.

\hs We introduce now  languages of infinite binary trees whose nodes
are labelled in a finite alphabet $\Si$.

\hs A node of an infinite binary tree is represented by a finite  word over
the alphabet $\{l, r\}$ where $r$ means ``right" and $l$ means ``left". Then an
infinite binary tree whose nodes are labelled  in $\Si$ is identified with a function
$t: \{l, r\}^\star \ra \Si$. The set of  infinite binary trees labelled in $\Si$ will be
denoted $T_\Si^\om$.
A tree language is a subset of $T_\Si^\om$, for some alphabet $\Si$.
(Notice that we shall only consider in the sequel {\it infinite} trees so we shall
often simply call  tree an {\it infinite} tree).

\hs Let $t$ be a tree. A branch $B$ of $t$ is a subset of the set of nodes of $t$ which
is linearly ordered by the tree partial order $\sqsubseteq$ and which
is closed under prefix relation,
i.e. if  $x$ and $y$ are nodes of $t$ such that $y\in B$ and $x \sqsubseteq y$ then $x\in B$.
\nl A branch $B$ of a tree is said to be maximal iff there is not any other branch of $t$
which strictly contains $B$.

\hs     Let $t$ be an infinite binary tree in $T_\Si^\om$. If $B$ is a maximal branch of $t$,
then this branch is infinite. Let $(u_i)_{i\geq 0}$ be the enumeration of the nodes in $B$
which is strictly increasing for the prefix order.
\nl  The infinite sequence of labels of the nodes of  such a maximal
branch $B$, i.e. $t(u_0)t(u_1) \ldots t(u_n) \ldots $  is called a path. It is an $\om$-word
over the alphabet $\Si$.

\hs We are now going to define tree automata and regular tree languages.

\begin{Deff} A (nondeterministic) tree automaton  is a quadruple $\mathcal{A}=(Q,\Si,\Delta, q_0)$, where $Q$
is a finite set of states, $\Sigma$ is a finite input alphabet, $q_0 \in Q$ is the initial state
and $\Delta \subseteq  Q \times   \Si   \times  Q \times   Q$ is the transition relation. 
 A run of the tree automaton  $\mathcal{A}$ on an infinite binary tree $t\in T_\Si^\om$ is an infinite binary tree $\rho \in T_Q^\om$ such that:
\nl (a)  $\rho (\varepsilon)=q_0$  and  ~~(b) for each $u \in \{l, r\}^\star$,  $(\rho(u), t(u), \rho(u.l), \rho(u.r))  \in \Delta$.
\end{Deff}

\begin{Deff}
A Muller  (nondeterministic) tree automaton  is a  5-tuple $\mathcal{A}=(Q,\Si,\Delta, q_0, \mathcal{F})$, where $(Q,\Si,\Delta, q_0)$ is a
tree automaton and $\mathcal{F} \subseteq 2^Q$ is the collection of  designated state sets.
 A run $\rho$ of the  Muller  tree automaton $\mathcal{A}$ on an infinite binary tree $t\in T_\Si^\om$ is said to be accepting if
for each path $p$ of $\rho$, the set of   states appearing infinitely  often on this path is in $\mathcal{F}$.
 The tree language $L(\mathcal{A})$ accepted by the  Muller tree automaton $\mathcal{A}$ is the set of infinite binary trees $t\in T_\Si^\om$
such that there is (at least) one accepting run of $\mathcal{A}$ on $t$.
 A tree language $L \subseteq T_\Si^\om$  is regular iff there exists a  Muller automaton $\mathcal{A}$ such that $L=L(\mathcal{A})$.
\end{Deff}

\begin{Rem}
A tree language is accepted by a Muller tree  automaton iff it is accepted by some Rabin tree automaton. We refer for instance to
\cite{Thomas90,PerrinPin} for the definition of Rabin tree  automaton.
\end{Rem}

\noi We now recall some fundamental closure properties of regular $\om$-languages and of regular tree languages.

\begin{The}[see \cite{Thomas90,PerrinPin}]
The class of regular  $\om$-languages (respectively, of regular tree languages) is effectively closed under finite union, finite intersection, and complementation, i.e.
we can effectively construct, from two B\"uchi automata  (respectively, Muller  tree automata)   $\mathcal{A}$ and $\mathcal{B}$, some
B\"uchi automata  (respectively, Muller tree automata)
$\mathcal{C}_1$, $\mathcal{C}_2$, and $\mathcal{C}_3$,  such that $L(\mathcal{C}_1)=L(\mathcal{A}) \cup L(\mathcal{B})$,
$L(\mathcal{C}_2)=L(\mathcal{A}) \cap L(\mathcal{B})$,  and $L(\mathcal{C}_3)$ is the complement of $L(\mathcal{A})$.
\end{The}

\noi Notice that one can  consider a relation $R \subseteq \Si_1^\om \times \Si_2^\om  \times \ldots \times \Si_n^\om$, where $\Si_1, \Si_2, \ldots \Si_n$, are  finite alphabets,
as an $\om$-language over the product alphabet $\Si_1 \times \Si_2  \times \ldots \times \Si_n$.  In a similar way we can consider a relation
$R \subseteq T_{\Si_1}^\om \times T_{\Si_2}^\om  \times \ldots \times T_{\Si_n}^\om$,
as a tree language over the product alphabet $\Si_1 \times \Si_2 \times \ldots \times \Si_n$.

\hs Let now $\mathcal{M}=(M, (R_i^M)_{1\leq i\leq n})$  be  a relational structure,
  where $M$ is the domain,  and for each $i\in [1, n]$ ~ $R_i^M$ is a relation
of finite arity $n_i$ on the domain $M$. The structure is said to be  $\om$-automatic (respectively,  $\om$-tree-automatic)
if there is a presentation of the structure
where the domain and the relations on the domain are accepted by B\"uchi automata (respectively, by Muller tree automata), in the following sense.

\begin{Deff}[see \cite{Blumensath99}]
Let $\mathcal{M}=(M, (R_i^M)_{1\leq i\leq n})$ be a relational structure, where $n\geq 1$ is an integer,  and each relation $R_i$ is of finite arity $n_i$.
\nl An  $\om$-tree-automatic presentation of the structure $\mathcal{M}$  is formed by   a tuple of  Muller tree  automata
$(\mathcal{A}, \mathcal{A}_=,  (\mathcal{A}_i)_{1\leq i\leq n})$,  and a mapping $h$ from $L(\mathcal{A})$ onto $M$,  such that:
\begin{enumerate}
\ite The automaton $\mathcal{A}_=$ accepts
an equivalence relation $E_\equiv $  on $L(\mathcal{A})$,  and
\ite
For each $i \in [1, n]$, the automaton $\mathcal{A}_i$ accepts an $n_i$-ary relation $R'_i$ on
$L(\mathcal{A})$ such that $E_\equiv$ is compatible with $R'_i$, and
\ite   The mapping $h$ is an isomorphism from the quotient  structure \nl $( L(\mathcal{A}),  (R'_i)_{1 \leq i \leq n} ) / E_\equiv$ onto $\mathcal{M}$.
\end{enumerate}

\noi  The $\om$-tree-automatic presentation is said to be injective if the equivalence relation $E_\equiv $ is just the equality relation on 
$L(\mathcal{A})$. In this case  $\mathcal{A}_=$ and $E_\equiv $  can be omitted and 
$h$ is simply an isomorphism from $( L(\mathcal{A}),  (R'_i)_{1 \leq i \leq n} )$ onto $\mathcal{M}$.
\noi A relational structure is said to be (injectively) $\om$-tree-automatic if it has an (injective) $\om$-tree-automatic presentation.
\end{Deff}

\noi Notice that sometimes an  $\om$-tree-automatic presentation is only given by a tuple of Muller tree  automata
$(\mathcal{A}, \mathcal{A}_=,  (\mathcal{A}_i)_{1\leq i\leq n})$, i.e. {\it without the mapping} $h$. In that case we still get the
$\om$-tree-automatic structure $( L(\mathcal{A}),  (R'_i)_{1 \leq i \leq n})  / E_\equiv$ which is in fact {\bf equal to } $\mathcal{M}$ {\bf up to isomorphism}.

\hs Notice also that, due to the good decidability properties of  Muller tree  automata, we can decide whether a given  automaton
$\mathcal{A}_=$ accepts
an equivalence relation $E_\equiv $  on $L(\mathcal{A})$ and whether, for each $i \in [1, n]$,  the
automaton $\mathcal{A}_i$ accepts an $n_i$-ary relation $R'_i$ on
$L(\mathcal{A})$ such that $E_\equiv$ is compatible with $R'_i$.

\hs We get the definition of $\om$-automatic (injective) presentation of a  structure and of  $\om$-automatic structure  by replacing simply
Muller tree  automata by B\"uchi automata in the above definition.

\hs We recall now two important properties of automatic structures.

\begin{The}[see \cite{Blumensath99}]
The class of  $\om$-tree-automatic  (respectively,  $\om$-automatic) structures is closed under first-order interpretations. In other words  if 
$\mathcal{M}$ is an $\om$-tree-automatic  (respectively,  $\om$-automatic) structure and $\mathcal{M}'$ is a relational structure which is
first-order interpretable in the structure $\mathcal{M}$, then the structure  $\mathcal{M}'$ is also $\om$-tree-automatic  (respectively,  $\om$-automatic).
\end{The}

\begin{The}[see \cite{Hodgson,Blumensath99}]\label{dec}
The first-order theory of an $\om$-tree-automatic (respectively,  $\om$-automatic) structure is decidable.
\end{The}

\noi Notice that $\om$-(tree)-automatic structures are always {\it relational} structures. However we can also consider structures equipped with
functional operations like groups, by replacing as usually a $n$-ary function by its graph which is a $(n+1)$-ary relation.
This will be always the case in the sequel where all structures are viewed as relational  structures.

\hs  Some examples of  $\om$-automatic structures can be found in
\cite{RubinPhd,NiesBSL,KNRS,KhoussainovR03,BlumensathGraedel04,KuskeLohrey,HjorthKMN08}.

\hs A first one is simply the boolean algebra of subsets of $\mathbb{N}$. The boolean algebra $\mathcal{P}(\mathbb{N})$ has an injective $\om$-automatic presentation
where any subset  $P\subseteq \mathbb{N}$ is simply represented by an infinite word $x_P$ over the alphabet $\{0,1\}$ defined by
$x_P(i)=1$ iff $i-1 \in P$ for all integers $i \in \mathbb{N}$. It is easy to see that the inclusion relation is then definable by a B\"uchi automaton, as well as the
operations of union, intersection, and complementation.

\hs The additive group $(\mathbb{R}, +)$ is $\om$-automatic, as is the product $(\mathbb{R}, +) \times (\mathbb{R}, +)$.

\hs Assume that a finite alphabet $\Si$ is linearly ordered.
Then the  set  $(\Sio, \leq_{lex})$ of infinite words over the alphabet $\Si$, equipped with the lexicographic ordering, is also $\om$-automatic.

\hs Is is easy to see that every (injectively) $\om$-automatic structure is also
(injectively) $\om$-tree-automatic. Indeed a Muller tree automaton can easily simulate a B\"uchi automaton on the leftmost branch of an infinite  tree.
The converse is not true, as shown in  \cite{HjorthKMN08} with   the following example.

\hs Let $\mathbb{T}=\{l, r\}^\star$. An element of $\mathbb{T}$ is simply a finite word over the alphabet 
$\{l, r\}$. A subset $U$ of $\mathbb{T}$ can be identified with an infinite binary
tree $t_U \in T_{\{0, 1\}}^\om$ such that for all $u\in \{l, r\}^\star$ it holds that $t_U(u)=1$ if and only if $u\in U$.
\nl Let then $V=\{ B \subseteq  \mathbb{T} \mid   t_B \notin  Path( (0^\star.1)^\om )  \}$, where we denote $Path((0^\star.1)^\om)$  the set of
infinite trees $t$ in $T_{\{0, 1\}}^\om$ such that $t$ has (at least) one path in $(0^\star.1)^\om$.
The structure $( \mathcal{P}(\mathbb{T}), \subseteq, V)$ is $\om$-tree-automatic but not  $\om$-automatic, see \cite{HjorthKMN08}.

\section{Two $\om$-tree automatic boolean algebras}

\noi We introduce in this section  two
$\om$-tree-automatic boolean algebras $\mathcal{B}_1$ and
$\mathcal{B}_2$ and we  will show later in Section  6 that the statement ``$\mathcal{B}_1$ is
isomorphic to $\mathcal{B}_2$" is independent from the axiomatic
system {\rm  ZFC}.

\hs Some examples of $\om$-automatic, hence also $\om$-tree-automatic, structures have been already given above.
Another example is the boolean algebra $\mathcal{P}(\mathbb{N})/Fin$ of subsets of $\mathbb{N}$ modulo {\it finite sets}.
The set $Fin$ of finite subsets of $\mathbb{N}$ is  an ideal of  $\mathcal{P}(\mathbb{N})$, i.e.  a subset of the powerset
of $\mathbb{N}$ such that:
\begin{enumerate}
\ite $\emptyset \in Fin$ and $\mathbb{N} \notin Fin$.
\ite For all $B, B' \in Fin$, it holds that $B \cup B' \in Fin$.
\ite For all $B, B' \in \mathcal{P}(\mathbb{N})$, if $B \subseteq B'$ and $B' \in Fin$ then $B \in Fin$.
\end{enumerate}

\noi
For any two subsets $A$ and $B$ of $\mathbb{N}$ we denote $A \Delta B$ their symmetric difference.
Then the relation $\approx$ defined by:
``$A \approx B$ iff the symmetric difference $A \Delta B$ is finite" is an equivalence relation on $\mathcal{P}(\mathbb{N})$.
The quotient
$\mathcal{P}(\mathbb{N}) / \approx $ denoted  $\mathcal{P}(\mathbb{N}) / Fin$ is a boolean algebra.
It is easy to see that this boolean algebra is $\om$-automatic.
We denote $[A]$ the equivalence class of a set of integers $A$. Let $\Si=\{0, 1\}$ and $L(\mathcal{A})= \Sio$ and for any $x\in \Sio$,
$h(x)=[\{ i\in \mathbb{N}\mid x(i+1)=1\}]$.
Then it is easy to see that $\{ (u, v) \in (\Sio)^2  \mid h(u)=h(v)  \}$ is accepted by a B\"uchi automaton.
Similarly the relation
$\{ (u, v) \in (\Sio)^2   \mid  h(u) \subseteq^\star h(v)  \}$ is a regular $\om$-language because the ``almost inclusion" relation $\subseteq^\star$
satisfies that $ h(u) \subseteq^\star h(v)$ iff $u(i) \leq v(i)$ for almost all integers $i$.

\hs The operations    $\cap, \cup, \neg$,
of intersection, union, and complementation,   on  $\mathcal{P}(\mathbb{N})/ Fin$ are defined by: 
  $[B]   \cap [B'] = [B  \cap  B'] $,    $[B]   \cup [B'] = [B  \cup  B'] $,  and
 $\neg [B] = [\neg B] $, see \cite{Jech}.

\hs The operations  of intersection, union, (respectively, complementation),
considered as ternary relations (respectively, binary relation)  are also
given by regular $\om$-languages.

\hs ${\bf 0}=[\emptyset]$ is the equivalence class of the empty set and ${\bf 1}=[\mathbb{N}]$ is the class of $\mathbb{N}$.

\hs The structures $(\mathcal{P}(\mathbb{N}) / Fin, \cap, \cup, \neg, {\bf 0},  {\bf 1})$ and
$(\mathcal{P}(\mathbb{N}) / Fin, \subseteq^\star)$ are $\om$-automatic, hence also
 $\om$-tree-automatic.

\hs We are going to consider now another boolean algebra.
 Let $\mathbb{T}=\{l, r\}^\star$ be the set of   finite words over the alphabet $\{l, r\}$. A subset $B$ of $\mathbb{T}$ has no infinite antichain (for the prefix order)
iff there is no infinite subset $D$ of $B$ such that for all $u, v \in D$, with $u\neq v$,
$u$ and $v$ are incomparable for the prefix order relation $\sqsubseteq$.

\hs
Let then $I=\{  B \subseteq \{l, r\}^\star \mid  B \mbox{ has no infinite antichain} \}$.  
The set $I$ is an ideal of  $\mathcal{P}(\mathbb{T})$, i.e. it is a subset of the powerset
of $\mathbb{T}$ such that:
\begin{enumerate}
\ite $\emptyset \in I$ and $\mathbb{T} \notin I$.
\ite For all $B, B' \in I$, it holds that $B \cup B' \in I$.
\ite For all $B, B' \in \mathcal{P}(\mathbb{T})$, if $B \subseteq B'$ and $B' \in I$ then $B \in I$.
\end{enumerate}

\noi  We can now consider the quotient $\mathcal{P}(\mathbb{T})/ I$ of the set of subsets of $\mathbb{T}$ modulo the ideal $I$. The relation $\approx_I$  defined on $\mathcal{P}(\mathbb{T})$ by:
``$A \approx_I B$ iff the symmetric difference $A \Delta B$ is in I" is an equivalence relation on $\mathcal{P}(\mathbb{T})$. The quotient
$\mathcal{P}(\mathbb{T}) / I$ is a boolean algebra.

\hs We are going to show that this boolean algebra is actually $\om$-tree-automatic.
Recall that a subset $B$ of $\mathbb{T}$ can be identified with an infinite binary
tree $t_B \in T_{\{0, 1\}}^\om$ such that for all $u\in \{l, r\}^\star$ it holds that $t_B(u)=1$ if and only if $u\in B$.
We can now state the following easy lemma.

\begin{Lem}\label{lemma-antichain}
Let $B \subseteq \{l, r\}^\star$. Then the set $B$ has an infnite antichain if and only if the tree $t_B$ has an  infinite branch, whose nodes form a sequence
$\varepsilon = u_0 \sqsubseteq u_1 \sqsubseteq u_2 \sqsubseteq \ldots$ and there exist infinitely many integers $n_i$ such that:
\nl   $u_{n_i}.a \sqsubseteq v_i$~~  for  some $a\in \{l, r\}$ and $v_i \in B$, and
\nl  $ u_{n_i +1} = u_{n_i}.b$~~  for $b\in \{l, r\}$ and $b\neq a$.
\end{Lem}

\proo Assume first that a tree $t_B$ has an  infinite branch, whose nodes form a sequence
$\varepsilon = u_0 \sqsubseteq u_1 \sqsubseteq u_2 \sqsubseteq \ldots$ and there exist infinitely many integers $n_i$ such that:
\nl   $u_{n_i}.a \sqsubseteq v_i$~~  for  some $a\in \{l, r\}$ and $v_i \in B$, and
\nl  $u_{n_i +1} = u_{n_i}.b$~~  for $b\in \{l, r\}$ and $b\neq a$.
\nl Then it is easy to see that the nodes $v_i \in B$ form an infinite antichain of $B$.

\hs Conversely, assume that $B \subseteq \{l, r\}^\star$ has an infinite antichain formed by nodes $w_i$, $i\geq 1$.
We can easily construct by induction an infinite sequence of nodes
$\varepsilon = u_0 \sqsubseteq u_1 \sqsubseteq u_2 \sqsubseteq \ldots$ forming an infinite branch $b$ and such that for each integer $i\geq 1$ there
are infinitely many nodes $w_j$ such that  $u_i \sqsubseteq w_j$.
Assume now, towards a contradiction, that there are only finitely many integers $n_i$ such that:
\nl   $u_{n_i}.a \sqsubseteq w_{l_i}$~~  for  some $a\in \{l, r\}$ and $l_i\geq 1$, and
\nl  $u_{n_i +1} = u_{n_i}.b$~~  for $b\in \{l, r\}$ and $b\neq a$.
\nl Then there exists an integer $N$ which is greater than all these integers $n_i$. Consider the node $u_N$ of the branch $b$. By construction
there are infinitely many integers $j$ such that $u_N \sqsubseteq w_j$. But all these nodes $w_j$ should be on the branch $b$. This would imply
that the nodes $w_i$, $i\geq 1$, do not form an infinite antichain. Thus this would lead to a contradiction so
there are  {\it infinitely many} integers $n_i$ such that:
\nl   $u_{n_i}.a \sqsubseteq w_{l_i}$~~  for  some $a\in \{l, r\}$ and $l_i\geq 1$, and
\nl  $u_{n_i +1} = u_{n_i}.b$~~  for $b\in \{l, r\}$ and $b\neq a$.
\ep

\hs  We can now state the following result.

\begin{Lem}
The set  $T=\{  t_B \mid  B \subseteq \{l, r\}^\star  \mbox{ and }  B \mbox{ has an  infinite antichain} \}$  is a regular tree language.
\end{Lem}

\proo We can construct a Muller tree automaton $\mathcal{A}$
accepting the set $T$. We explain informally the behaviour of this automaton. Using the
non-determinism the automaton, when reading a tree $t_B$, will guess an infnite branch whose existence is given by the preceding lemma, and integers $n_i$
having the property given by the same lemma.
\ep

\hs We know that  the class of regular tree languages
 is effectively closed under complementation. So we get now the following result.

\begin{Lem}\label{T_I}
The set  $T_I=\{  t_B \mid  B \subseteq \{l, r\}^\star  \mbox{ and }  B \mbox{ has no  infinite antichain} \}$  is a regular tree language.
\end{Lem}

\noi 
 Recall that, if we denote by $[A]$ the equivalence class of a set $A \subseteq \mathbb{T}=\{l, r\}^\star$, then  the operations
of intersection, union, and complementation,   on  $\mathcal{P}(\mathbb{T})/ I$ 
are defined by:   $[B]   \cap [B'] = [B  \cap  B'] $,    $[B]   \cup [B'] = [B  \cup  B'] $,  and
 $\neg [B] = [\neg B] $. The almost inclusion relation $\subseteq^\star$ is defined by:  $[B]   \subseteq^\star [B']$ iff  $B \setminus B'  \in I$.

\hs We can now state the following result which will be fundamental in the sequel.

\begin{Pro}
The boolean algebra $(\mathcal{P}(\mathbb{T})/ I,  \cap, \cup, \neg, {\bf 0},  {\bf 1})$ and the structure $(\mathcal{P}(\mathbb{T})/ I, \subseteq^\star)$  are $\om$-tree-automatic structures.
\end{Pro}

\proo
We denote by $[A]$ the equivalence class of a set $A \subseteq \mathbb{T}=\{l, r\}^\star$.
The class $[A]$ will be represented by the trees $t_B   \in T_{\{0, 1\}}^\om$ such that $[A]=[B]$, i.e. such that the symmetric difference $A \Delta B$ is in I.
The function $h$ will associate the class   $[B]$     to each tree $t_B \in T_{\{0, 1\}}^\om$.
Then it is easy to see, from Lemma \ref{T_I},
 that $\{ (t, t') \in T_{\{0, 1\}}^\om \times T_{\{0, 1\}}^\om  \mid h(t)=h(t')  \}$ is accepted by a Muller tree automaton accepting
infinite trees in $T_{\{0, 1\}\times \{0, 1\}}^\om$.
\nl From Lemma \ref{T_I}, we can also easily infer that the almost inclusion relation and the operations
of intersection, union, and complementation,    are also
given by regular tree languages.
${\bf 0}=[\emptyset]$ is simply represented by the trees in  the set $T_I$. The class ${\bf 1}=[\mathbb{T}]$ is represented by the trees $t_B$, where
$B \subseteq \{l, r\}^\star$   and  $\neg B$ has no  infinite antichain.
\ep

\hs From now on we shall  denote  $\mathcal{B}_1 = (\mathcal{P}(\mathbb{N}) / Fin, \cap, \cup, \neg, {\bf 0},  {\bf 1})$
and $\mathcal{B}_2 = (\mathcal{P}(\mathbb{T})/ I,  \cap, \cup, \neg, {\bf 0},  {\bf 1})$  the two $\om$-tree-automatic boolean algebras defined above.

\hs Recall now the definition of an atomless boolean algebra. 

\begin{Deff}
Let $\mathcal{B} = (B, \cap, \cup, \neg, {\bf 0},  {\bf 1})$ be a boolean algebra and $\subseteq$ be the inclusion relation on $B$ defined by 
$x \subseteq y$ iff $x \cap y = x$ for all $x, y \in B$. Then the boolean algebra $\mathcal{B}$ is said to be an atomless boolean algebra iff for every 
$x\in B$ such that $x\neq {\bf 0}$ there exists a $z\in B$ such that ${\bf 0} \subset z \subset x$. 
\end{Deff}

\noi We can now state the following result. 

\begin{Pro} The two boolean algebras $\mathcal{B}_1$ and $\mathcal{B}_2 $ are atomless boolean algebras. 
\end{Pro}

\proo  Consider firstly the  boolean algebra $\mathcal{B}_1 = (\mathcal{P}(\mathbb{N}) / Fin, \cap, \cup, \neg, {\bf 0},  {\bf 1})$. 
It is well known that it is an atomless boolean algebra. We now  give a proof of  this result. 
Let $A \subseteq \mathbb{N}$ be  such that the equivalence class $[A]$ is different from
the element ${\bf 0}$ in $\mathcal{B}_1$. Then the set $A$ is infinite and there exist two infinite sets $A_1$ and $A_2$ such that $A=A_1 \cup A_2$. 
The element $[A_1]$ is different from the element ${\bf 0}$ in $\mathcal{B}_1$ because $A_1$ is infinite, and $[A_1] \subset [A]$ because 
$A - A_1=A_2 $ is infinite. Thus the following strict inclusions hold in   $\mathcal{B}_1$: 
${\bf 0} \subset [A_1]  \subset [A]$. 

\hs   We can prove in a similar way that the boolean algebra $\mathcal{B}_2 $ is an  atomless boolean algebra. Let then 
$X \subseteq \mathbb{T}$  be  such that the equivalence class $[X]$ is different from
the element ${\bf 0}$ in $\mathcal{B}_2$. Then the set $X \subseteq \{l, r\}^\star$ contains  an infinite antichain $A\subseteq X$.  
There are two infinite sets $A_1$ and $A_2$ such that $A=A_1 \cup A_2$. These two sets are also infinite antichains, so $[A_1]$ is different from 
${\bf 0}$ in $\mathcal{B}_2$ and $[A_1] \subset [X]$ because the set $X - A_1$ contains the infinite antichain $A_2$. 
Thus the following strict inclusions hold in   $\mathcal{B}_2$: 
${\bf 0} \subset [A_1]  \subset [X]$. This proves that the boolean algebra $\mathcal{B}_2$ is atomless. 
\ep 

\section{Topology}

\noi We assume the reader to be familiar with basic notions of topology which
may be found in \cite{Moschovakis80,LescowThomas,Kechris94,Staiger97,PerrinPin}.
There is a natural metric on the set $\Sio$ of  infinite words 
over a finite alphabet 
$\Si$ containing at least two letters which is called the {\it prefix metric} and defined as follows. For $u, v \in \Sio$ and 
$u\neq v$ let $\delta(u, v)=2^{-l_{\mathrm{pref}(u,v)}}$ where $l_{\mathrm{pref}(u,v)}$ 
 is the first integer $n$
such that $u(n+1)\neq v(n+1) $. 
This metric induces on $\Sio$ the usual  Cantor topology for which {\it open subsets} of 
$\Sio$ are in the form $W.\Si^\om$, where $W\subseteq \Sis$.
A set $L\subseteq \Si^\om$ is a {\it closed set} iff its complement $\Si^\om - L$ 
is an open set.  We recall now a characterization of closed sets which will be useful in the sequel. 

\begin{Pro}[see \cite{Staiger97,Kechris94}]\label{carac-closed}
A set $L\subseteq \Si^\om$ is a closed subset of $\Si^\om$ iff for every $\sigma\in \Si^\om$, 

$[\fa n\geq 1,  \exists u\in \Si^\om$  such that $\sigma (1)\ldots \sigma (n).u \in L]$
 implies that $\sigma\in L$.
\end{Pro}

\noi We  now define the next classes of the Borel Hierarchy of  subsets of $\Sio$.

\begin{Deff}
For an integer $n\geq 1$,  the classes ${\bf \Si}^0_n$
 and ${\bf \Pi}^0_n$ of the Borel Hierarchy on the topological space $\Si^\om$ 
are defined as follows:
\nl ${\bf \Si}^0_1$ is the class of open subsets of $\Si^\om$, 
 ${\bf \Pi}^0_1$ is the class of closed subsets of $\Si^\om$, 
\nl and for any  integer $n\geq 1$: 
\nl ${\bf \Si}^0_{n+1}$ is the class of countable unions of ${\bf \Pi}^0_n$-subsets of $\Si^\om$. 
 \nl ${\bf \Pi}^0_{n+1}$ is the class of countable intersections of ${\bf \Si}^0_n$-subsets of $\Si^\om$.
\end{Deff}

\begin{Notation} 
Following the earlier notations for the Borel hierarchy of Borel sets of finite rank,  a ${\bf \Pi}^0_2 $-set is also called a    $G_\delta$-set 
and a  ${\bf \Si}^0_2$-set is also called a  $F_\sigma$-set.  So a $G_\delta$-set is a countable intersection of open sets and a $F_\sigma$-set 
is a countable union of closed sets. 
\end{Notation}

\noi   The Borel Hierarchy is also defined for transfinite levels indexed by countable ordinals, see \cite{Moschovakis80,Kechris94}. 
 However we shall not need these notions in the sequel. Recall that the class of Borel subsets of $\Si^\om$ is the closure of the class of open subsets of $\Si^\om$ 
under countable union and countable intersection. 

\hs  
There exists another hierarchy beyond the Borel hierarchy, which is called the 
projective hierarchy. 
The first level of the projective hierarchy is formed by the class ${\bf \Si}^1_1$  of  {\it analytic sets} and the 
class ${\bf \Pi}^1_1$ of {\it co-analytic sets} which are complements of 
analytic sets.  
In particular 
the class of Borel subsets of $\Si^\om$ is strictly included in  
the class  ${\bf \Si}^1_1$ of {\it analytic sets} which are 
obtained by projection of Borel sets. 

\begin{Deff} 
A subset $A$ of  $\Si^\om$ is in the class ${\bf \Si}^1_1$ of {\bf analytic} sets
iff there exists a finite set $Y$ and a Borel subset $B$  of  $(\Si \times Y)^\om$ 
such that $ x \in A \lra [\exists y \in Y^\om ~~ (x, y) \in B]$, 
where $(x, y)$ is the infinite word over the alphabet $\Si \times Y$ such that
$(x, y)(i)=(x(i),y(i))$ for each  integer $i\geq 1$.
\end{Deff} 

\noi  An important fact in this paper is   that  the powerset $\mathcal{P}(\mathbb{N})$ can be  equipped with the 
standard metric topology obained from its identification with the Cantor space $\{0, 1\}^\om$.  
Then the  topological notions like open, closed, $F_\sigma$, analytic, can be applied 
to families of subsets of $\mathbb{N}$. 
 
\hs The ideal $Fin$ of  $\mathcal{P}(\mathbb{N})$ is  identified with the set of $\om$-words over the alphabet $\{0, 1\}$ having only finitely many 
letters $1$. It is a well known example of  $F_\sigma$-subset of $\{0, 1\}^\om$, as stated in the following lemma. 

\begin{Lem}\label{fin}[see \cite{PerrinPin,Kechris94}]
The powerset $\mathcal{P}(\mathbb{N})$ being identified with the Cantor space $\{0, 1\}^\om$, the ideal $Fin$ of  $\mathcal{P}(\mathbb{N})$ 
is a $F_\sigma$-subset of $\{0, 1\}^\om$. 
\end{Lem}

\proo Let $k\geq 1$ be an integer and $Fin_k$ be the set of subsets of $\mathbb{N}$ having at most $k$ elements. The set 
$Fin_k$ is identified with the set of $\om$-words over the alphabet $\{0, 1\}$ having at most $k$ letters $1$. Using Proposition \ref{carac-closed} it 
is easy to see that for each $k\geq 1$ the set $Fin_k$ is then a closed subset of $\{0, 1\}^\om$. This follows from the fact that if a set 
$A \subseteq \mathbb{N}$ is such that every finite subset of $A$ has at most $k$ elements, then the set $A$ itself  has at most $k$ elements.
Thus $Fin=\bigcup_{k\geq 1} Fin_k$ 
is a countable union of closed sets,  i.e. a $F_\sigma$-subset of $\{0, 1\}^\om$.
\ep

\hs  Consider now the set $\mathbb{T}=\{l, r\}^\star$ of  finite words over the alphabet $\{l, r\}$. This set is countably infinite and we can define a bijection from 
$\mathbb{T}=\{l, r\}^\star$ onto $\mathbb{N}$ by enumerating the elements of $\mathbb{T}$.  
For each integer $n\geq 0$, call $W_n$ the set of words of length $n$ of $\{l, r\}^\star$. 
Then $W_0=\{\varepsilon\}$, $W_1=\{l, r\}$, $W_2=\{ll, lr, rl, rr\}$ and so on.
$W_n$ is the set of nodes which appear in the $(n+1)^{th}$ level of an infinite binary tree.
 We consider now  
the lexicographic order on $W_n$ (assuming that $l$ is before $r$ for this order).
Then, in the enumeration of the nodes with regard to this order, the  nodes of $W_1$ will 
be: $l, r$; the nodes of $W_3$ will be: $ lll, llr, lrl, lrr, rll, rlr, rrl, rrr$.
We enumerate now the  elements of $\mathbb{T}$ in the following order. We begin with $\varepsilon$, then the nodes in $W_1$ in the 
lexicographic order, then the nodes in $W_2$ in the 
lexicographic order, then  the nodes in $W_3$ in the 
lexicographic order, and so on \ldots
The successive nodes are then 

$$\varepsilon, l, r, ll, lr, rl, rr,  lll, llr, lrl, lrr, rll, rlr, rrl, rrr, \ldots$$

\noi For every $u \in \{l, r\}^\star$, we define $f(u) \in \mathbb{N}$ such that $u$ is the $(f(u)+1)th$  element in the above enumeration of words of 
$\{l, r\}^\star$. For instance $f(\varepsilon)=0$, $f(l)=1$, $f(r)=2$, $f(ll)=3$, $f(lr)=4$, and so on \ldots
\nl The function $f$ is then a bijection from $ \{l, r\}^\star$ onto $\mathbb{N}$, and it induces also a  bijection  from $\mathcal{P}(\mathbb{T})$ onto 
$\mathcal{P}( \mathbb{N})$ which will be also denoted by $f$, the meaning being clear from the context. 

\hs Recall that we have set $I=\{  B \subseteq \{l, r\}^\star \mid  B \mbox{ has no infinite antichain} \}$, and that the set $I$ is an ideal of  $\mathcal{P}(\mathbb{T})$. 
We are going to show that $f(I)$ is a $F_\sigma$-subset of $\{0, 1\}^\om$, where again  the powerset $\mathcal{P}(\mathbb{N})$  
is identified with the Cantor space $\{0, 1\}^\om$.  
We first state the following lemma. 

\begin{Lem}\label{union}
 Let $I=\{  B \subseteq \{l, r\}^\star \mid  B \mbox{ has no infinite antichain} \}$ and, for each 
integer  $k\geq 1$,  
$I_k=\{  B \subseteq \{l, r\}^\star \mid  B \mbox{ has no antichain of cardinal  greater}$ $\mbox{than}$ $k \}$.  
Then  $I=\bigcup_{k\geq 1} I_k$. 
\end{Lem}

\proo It is clear that if $B \subseteq \{l, r\}^\star$ has no antichain of cardinal  greater than $k$, for an integer $k\geq 1$, then $B$ has no infinite antichain. 
Then the inclusion $\bigcup_{k\geq 1} I_k \subseteq I$ holds. 

\hs  We want now to  prove that $I \subseteq \bigcup_{k\geq 1} I_k $.  
We assume that a set $B \subseteq \{l, r\}^\star$ is not in $\bigcup_{k\geq 1} I_k $ and 
we are going to  prove that this implies that  $B \notin I$. 

\hs Let then $B \subseteq \{l, r\}^\star$ such that $B \notin \bigcup_{k\geq 1} I_k$. For every $k\geq 1$ it holds that $B \notin I_k$, i.e.
 $B$ has some antichain of cardinal  greater than $k$. We are going to prove that $B$ has an infinite antichain, using the characterization  given by 
Lemma \ref{lemma-antichain}.  We can first construct by induction an infinite sequence of nodes of the tree $t_B$: 
$$\varepsilon = u_0 \sqsubseteq u_1 \sqsubseteq u_2 \sqsubseteq \ldots$$
\noi such that for every integer $j\geq 0$  and every  $k\geq 1$,  there is an antichain in $B$ of cardinal  greater than $k$ whose nodes have 
$u_j$ as prefix.
Indeed assume that we have already construct 
$\varepsilon = u_0 \sqsubseteq u_1 \sqsubseteq u_2 \sqsubseteq \ldots  \sqsubseteq u_j$. 
Then at least one node among $u_j.l$ and $u_j.r$ has the desired property 
and we can choose this node as $u_{j+1}$. 

\hs We can now see that there exist   infinitely many integers $n_i$ such that:
\nl   $u_{n_i}.a \sqsubseteq v_i$~~  for  some $a\in \{l, r\}$ and $v_i \in B$, and
\nl  $ u_{n_i +1} = u_{n_i}.b$~~  for $b\in \{l, r\}$ and $b\neq a$.

\hs Indeed we can construct the sequence $(n_i)_{i\geq 0}$ by induction. Firstly there is an antichain in $B$ of cardinal  greater than $2$ whose nodes have 
$u_0$ as prefix. Thus there is an integer $n_0\geq 0$ such that: 
$u_{n_0}.a \sqsubseteq v_0$~~  for  some $a\in \{l, r\}$ and $v_0 \in B$, and
  $ u_{n_0 +1} = u_{n_0}.b$~~  for $b\in \{l, r\}$ and $b\neq a$.
 Assume now that we have constructed integers $n_0, n_1, \ldots , n_j$ having the desired property. 
Then by construction of the sequence of nodes $(u_i)_{i\geq 0}$, we know that there is an antichain in $B$ of cardinal  greater than $2$ whose nodes have 
$u_{n_j +1}$ as prefix. Thus there is an integer $n_{j+1}\geq n_j +1$ such that: 
$u_{n_{j+1}}.a \sqsubseteq v_{j+1}$~~  for  some $a\in \{l, r\}$ and $v_{j+1}\in B$, and
  $ u_{n_{j+1} +1} = u_{n_{j+1}}.b$~~  for $b\in \{l, r\}$ and $b\neq a$.

\hs Using Lemma \ref{lemma-antichain} we can conclude that $B$ has an infinite antichain, i.e. $B \notin I$. This proves the 
inclusion $I \subseteq \bigcup_{k\geq 1} I_k $. 
\ep

\begin{Lem}\label{closed}
For each integer $k\geq 1$,  the set $f(I_k)$ is a closed subset of the Cantor space $\{0, 1\}^\om$.  
\end{Lem}

\proo Let $k\geq 1$ be an integer. We shall prove that $f(I_k)$ is a closed subset of the Cantor space $\{0, 1\}^\om$, using the 
characterization of closed sets given by Proposition \ref{carac-closed}. 

\hs Let $x \in \{0, 1\}^\om$ such that $[\fa n\geq 1,  \exists y_n \in \{0, 1\}^\om$  such that $x(1)\ldots x(n).y_n \in f(I_k)]$.  This implies that 
if $P_n$ is the subset of  $\mathbb{N}$ which is  identified with the $\om$-word $x(1)\ldots x(n).y_n$ then $f^{-1}(P_n)\subseteq \{l, r\}^\star$ 
 has no antichain of cardinal greater than $k$. In particular, for every $n\geq 1$, the finite set  $f^{-1}(P_n \cap \{0, \ldots, n-1\})$ 
 has no antichain of cardinal greater than $k$. Thus every finite subset of  $f^{-1}(x)$ has no antichain of cardinal greater than $k$. 
This implies that $f^{-1}(x)$ itself has no antichain of cardinal greater than $k$, i.e. $f^{-1}(x)$ is in $I_k$ and $x$ belongs to $f(I_k)$. 

\hs Using Proposition \ref{carac-closed} we can conclude that $f(I_k)$ is a closed subset of the Cantor space $\{0, 1\}^\om$.  
\ep

\hs We can now state the following result, which will be important in the sequel. 

\begin{Pro}\label{I}
Let  $I=\{  B \subseteq \{l, r\}^\star \mid  B \mbox{ has no infinite antichain} \}$. Then the set $f(I)$ is a $F_\sigma$-subset of $\{0, 1\}^\om$. 
\end{Pro}
 
\proo By Lemma \ref{union}, it holds that $I=\bigcup_{k\geq 1} I_k$. But by Lemma \ref{closed} for each  integer $k\geq 1$,  
the set $f(I_k)$ is a closed subset of the Cantor space $\{0, 1\}^\om$.  Thus $f(I)=\bigcup_{k\geq 1} f(I_k)$ is a countable 
union of closed sets, i.e. a $F_\sigma$-subset of $\{0, 1\}^\om$.
\ep 

\section{Axioms of set theory}

\noi We now recall some basic notions of set theory 
which will be useful in the sequel, and which are exposed in a  textbook on set theory, like \cite{Jech}. 

\hs  The usual axiomatic system ${\rm ZFC}$ is 
Zermelo-Fraenkel system  ${\rm ZF}$ plus the axiom of choice ${\rm AC}$. 
 A model ({\bf V}, $\in)$ of  the axiomatic system    ${\rm ZFC}$       is a collection  {\bf V} of sets,  equipped with 
the membership relation $\in$, where ``$x \in y$" means that the set $x$ is an element of the set $y$, which satisfies the axioms of  ${\rm ZFC}$.  
We shall often say `` the model {\bf V}"
instead of  ``the model  ({\bf V}, $\in)$".

\hs We recall that the infinite cardinals are usually denoted by
$\aleph_0, \aleph_1, \aleph_2, \ldots , \aleph_\alpha, \ldots$
\nl  and  that  Cantor's   Continuum
Hypothesis         ${\rm CH}$ states that the cardinality of the continuum $2^{\aleph_0}$ is equal
to the first uncountable cardinal $\aleph_1$. 

\hs Recall also that  ${\rm OCA}$ denotes the Open Coloring Axiom, a natural
alternative to ${\rm CH}$ which has been first considered by
the second author in \cite{Todorcevic89}. 

\hs
 For any set $X$ we denote   $[X]^2$ the set of subsets of $X$ having exactly $2$ elements. 
Let $\mathbb{R}$ be the set of reals equipped with the usual topology.  If $X   \subseteq \mathbb{R}$ and   $K  \subseteq   [X]^2$ 
then we say that $K$ is open if $\{ (x, y) \mid \{x, y\} \in K \}$ is open in $X \times X$. 

\hs The ${\rm OCA}$ states that if $X   \subseteq \mathbb{R}$  and $[X]^2 = K_0 \cup K_1$ is a partition of  $[X]^2$ with $K_0$ open, then either 
there exists an uncountable subset $Y$ of $X$ such that $[Y]^2 \subseteq K_0$ or there exist a sequence of sets $(H_n)_{n\in \om}$, 
such that $X=\bigcup_{n\in\om}H_n$ and,  for all $n \in \om$,  $[H_n]^2 \subseteq K_1$. 

\hs  The above version of the axiom ${\rm OCA}$  can be shown to be equivalent to the following  one: 
\nl If $G=(V, E)$ is a graph whose edge relation $E$ can be
written as a countable union of 'rectangles', i.e., sets of the
form
$$\{ \{x,y\}: x\in P, y\in Q \}$$ for $P,Q \subseteq V,$ then either
the chromatic number of $G$ is countable, or else $G$ has an
uncountable clique.  

\hs It is known that if the theory {\rm ZFC} is
consistent, then so are the theories {\rm (ZFC + CH)} and {\rm
(ZFC + OCA)}, see \cite[pages 176 and  577]{Jech}.

\hs We now recall some known properties of the  class ${\bf L}$ of  {\it constructible sets} in a model {\bf V} of ${\rm ZF}$, 
which will be useful in the sequel. 
If  {\bf V} is  a model of ${\rm ZF}$ and ${\bf L}$ is  the class of  {\it constructible sets} of   {\bf V}, then the class  ${\bf L}$     forms a model of  
${\rm (ZFC + CH)}$.
Notice that the axiom ({\bf V=L}) means ``every set is constructible"  and that it is consistent with ${\rm ZFC}$. 

\hs In particular, if {\bf V} is  a model of ${\rm (ZFC + OCA)}$ and if  ${\bf L}$ is  the class of  {\it constructible sets} of   {\bf V},  
then the class  ${\bf L}$     forms a model of  
${\rm (ZFC + CH)}$.

\section{The isomorphism relation}

\noi We are going to see that the statement ``$\mathcal{B}_1$ is
isomorphic to $\mathcal{B}_2$" is independent from the axiomatic
system {\rm  ZFC}.

\begin{The}\label{ind}
\noi
\begin{enumerate}
\ite  {\rm  (ZFC + CH)}   ~~$\mathcal{B}_1$ is isomorphic to
$\mathcal{B}_2$. 
\ite  {\rm   (ZFC + OCA) } $\mathcal{B}_1$ is not
isomorphic to $\mathcal{B}_2$.
\end{enumerate}
\end{The}

\proo 
We  have already seen  that the powerset $\mathcal{P}(\mathbb{N})$ can be  equipped with the 
standard metric topology obained from its identification with the Cantor space $\{0, 1\}^\om$.

\hs 
Then the ideal $Fin$ of  $\mathcal{P}(\mathbb{N})$ is  identified with the set of $\om$-words over the alphabet $\{0, 1\}$ having only finitely many 
letters $1$ and  Lemma \ref{fin} states that  it is a   $F_\sigma$-subset of $\{0, 1\}^\om$. Thus the boolean algebra $\mathcal{B}_1$ is a quotient algebra of
$\mathcal{P}(\mathbb{N})$ over a $F_\sigma$-ideal.  

\hs  Consider now the boolean algebra $\mathcal{B}_2= \mathcal{P}(\mathbb{T})/ I$. Recall that we have defined a bijection $f$ from 
$\mathbb{T}=\{l, r\}^\star$ onto $\mathbb{N}$ by enumerating the elements of $\mathbb{T}$.  Then the boolean algebra $\mathcal{B}_2$ is isomorphic to the 
boolean algebra $\mathcal{P}(\mathbb{N})/f(I)$. But by Proposition \ref{I}, the set $\mathcal{P}(\mathbb{N})$ being again 
identified with the Cantor set  $\{0, 1\}^\om$,  the set  $f(I)$ is a $F_\sigma$-subset of $\{0, 1\}^\om$.  Thus the boolean algebra $\mathcal{B}_2$ is also isomorphic to a quotient algebra of
$\mathcal{P}(\mathbb{N})$ over a $F_\sigma$-ideal.

\hs On the other hand,  Just and Krawczyk     proved in  \cite[Theorem 1]{JK} that two   boolean algebras which are quotients of 
$\mathcal{P}(\mathbb{N})$ over  $F_\sigma$-ideals are always isomorphic under  {\rm  (ZFC + CH)}. This implies the first part of the Theorem.

\hs On the other hand, it is proved in \cite[Theorem
6]{todorcevic98}(see also \cite{Farah} and \cite{Just} for more
information about the influence of ${\rm OCA}$ to the quotient
structures of this sort) that under {\rm OCA} if a quotient
algebra $\mathcal{P}(\mathbb{N})/J$ over an analytic ideal $J$ on
$\mathbb{N}$ is isomorphic to a subalgebra of
$\mathcal{P}(\mathbb{N})/Fin$ then $J$ must be a \emph{trivial
modification} of the ideal $Fin$, i.e., there is some infinite
subset $B$ of $\mathbb{N}$ such that
$$J=\{A \subseteq \mathbb{N}: A \cap B \in Fin\}.$$

\hs Consider  now  $I=\{  B \subseteq \{l, r\}^\star \mid  B \mbox{ has no infinite antichain}\}$ and  the ideal $J=f(I)$. By Proposition \ref{I} the ideal 
$f(I)$ is a $F_\sigma$ hence also analytic  subset  of $\{0, 1\}^\om$.  

\hs Let us now prove that  the ideal  $J=f(I)$    is not a trivial
modification of $Fin$. Towards a contradiction, assume on the contrary that there exists a  subset $B$ of $\mathbb{N}$ such that
$J=\{A \subseteq \mathbb{N}: A \cap B \in Fin\}.$ Notice that if $C \subseteq \{l, r\}^\star$ is a (finite or infinite) chain, i.e. is linearly ordered by 
the prefix order relation, then it has no infinite antichain and it belongs to $I$. In particular, for every integer $n\geq 1$ the set $C_n=\{l^n.r^k \mid k\geq 1\}$ 
is in $I$ so $f(C_n)$ is in $J=f(I)$ and $f(C_n)\cap B$ would be finite. This would imply that for every integer $n\geq 1$ there would exist an integer $k_n\geq 1$ such that 
$f(\{l^n.r^k \mid k\geq k_n\}) \cap B$ is empty. Let now $A=\{l^n.r^{k_n} \mid n\geq 1\}$. It is an infinite antichain so it is not in $I$. Thus $f(A)$ is not in $J$ 
 and $f(A)\cap B$ should be infinite. But by construction $f(A)\cap B$ is empty and this leads to a contradiction. This proves that the ideal  $J=f(I)$    is not a trivial
modification of $Fin$.

\hs 
Thus, assuming ${\rm OCA}$,  the boolean
algebra $\mathcal{P}(\mathbb{N})/f(I)$ is not even isomorphic to a subalgebra of
$\mathcal{B}_1.$ But the boolean algebra 
$\mathcal{B}_2$ is isomorphic to the boolean algebra $\mathcal{P}(\mathbb{N})/f(I)$ therefore 
the boolean
algebra $\mathcal{B}_2$  is not even isomorphic to a subalgebra of
$\mathcal{B}_1.$ 
\ep

\hs  We can now state the following result.

\begin{Cor}
The isomorphism relation for $\om$-tree-automatic structures
(respectively, $\om$-tree-automatic boolean algebras,
$\om$-tree-automatic partial orders) is not determined by the
axiomatic system  {\rm  ZFC}.
\end{Cor}

\proo The result for $\om$-tree-automatic boolean algebras,  hence also for $\om$-tree-automatic structures, follows directly from Theorem \ref{ind} and the fact
that the boolean algebras $\mathcal{B}_1$ and $\mathcal{B}_2$ are $\om$-tree-automatic.
For partial orders, we consider the $\om$-tree-automatic structures $(\mathcal{P}(\mathbb{N}) / Fin, \subseteq^\star)$  and $(\mathcal{P}(\mathbb{T})/ I, \subseteq^\star)$.
These two structures are isomorphic if and only if the two boolean algebras $\mathcal{B}_1$ and $\mathcal{B}_2$ are isomorphic, see \cite[page 79]{Jech}.
Then the result for partial orders  follows from the case of boolean algebras.
\ep

\hs We are going to get similar results for other classes of $\om$-tree-automatic structures. First we can consider a boolean algebra
$(B, \cap, \cup, \neg, {\bf 0},  {\bf 1})$
as a commutative ring with unity element    $(B, \Delta, \cap, {\bf 1})$, where $\Delta$ is the symmetric difference operation.
It is clear that the operations of union and complementation
can be defined from the symmetric difference and intersection operations.
Moreover two boolean algebras $(B, \cap, \cup, \neg, {\bf 0},  {\bf 1})$ and $(B', \cap, \cup, \neg, {\bf 0},  {\bf 1})$ are isomorphic if and only if
the rings $(B, \Delta, \cap, {\bf 1})$ and $(B', \Delta, \cap, {\bf 1})$ are isomorphic.
We shall denote $\mathcal{R}_1=(\mathcal{P}(\mathbb{N}) / Fin, \Delta, \cap,  {\bf 1})$  and
$\mathcal{R}_2=(\mathcal{P}(\mathbb{T})/ I, \Delta, \cap,  {\bf 1})$ the two commutative rings associated with the two boolean algebras
$\mathcal{B}_1$ and $\mathcal{B}_2$.
Notice that $[\emptyset]$ is the unity element for the operation $\Delta$ and that every element of the ring $\mathcal{R}_1$
(respectively, $\mathcal{R}_2$) is its own inverse for this group operation. On the other hand ${\bf 1}$ is the  unity element for the operation
$\cap$ in both rings and it is also the unique invertible element for the operation $\cap$ in both rings.

\hs
We can now state the following result, which follows directly from Theorem \ref{ind}.

\begin{The}\label{ind2}
\noi
\begin{enumerate}
\ite  {\rm  (ZFC + CH)}   ~$\mathcal{R}_1$ is isomorphic to
$\mathcal{R}_2$. \ite  {\rm   (ZFC + OCA) } $\mathcal{R}_1$ is not
isomorphic to $\mathcal{R}_2$.
\end{enumerate}
\end{The}

\noi  We can also obtain a similar result for non commutative rings. For that purpose we consider the set $M_n(R)$ of square matrices with $n$ columns and $n$
rows and coefficients in a given ring $R$.  If $n\geq 2$ then the set $M_n(R)$, equipped with addition and multiplication of matrices, is a non commutative ring.
The ring $M_n(R)$ is first-order interpretable in the ring $R$; each matrix $M$ being represented by a unique $n^2$-tuple of elements of $R$,  the addition and
multiplication of matrices are first order definable in $R$.

\hs On the other hand,  the class of  $\om$-tree-automatic structures is closed under first order interpretations.
Thus if $R$ is an  $\om$-tree-automatic ring   then the ring of matrices
$M_n(R)$ is also $\om$-tree-automatic.
We now denote $\mathcal{M}_1=M_n(\mathcal{R}_1)$ and $\mathcal{M}_2=M_n(\mathcal{R}_2)$, for some fixed integer $n\geq 2$,
   and state the following result.

\begin{The}\label{ind3}
\noi
\begin{enumerate}
\ite  {\rm  (ZFC + CH)}   ~~$\mathcal{M}_1$ is isomorphic to
$\mathcal{M}_2$. \ite  {\rm   (ZFC + OCA) } $\mathcal{M}_1$ is not
isomorphic to $\mathcal{M}_2$.
\end{enumerate}
\end{The}

\proo  It is clear that if $\mathcal{R}_1$ is isomorphic to $\mathcal{R}_2$ then the rings $\mathcal{M}_1$ and  $\mathcal{M}_2$ are also isomorphic.
Conversely, assume that $\Phi  :  \mathcal{M}_1 \ra \mathcal{M}_2$ is an  isomorphism of rings. Then, for $i \in \{1, 2\}$,
consider the center $\mathcal{C}_i$ of the ring   $\mathcal{M}_i$. The set $\mathcal{C}_i$ contains the matrices of $\mathcal{M}_i$
which commute with every matrix of $\mathcal{M}_i$.  It holds that the restriction of  $\Phi$ to $\mathcal{C}_1$ is an isomorphism between $\mathcal{C}_1$
and $\mathcal{C}_2$. But it is well known that the center $C$ of a ring $M_n(R)$  is formed by matrices $M(u)$
which have a same element $u\in R$ on the diagonal and zeros elsewhere.
Then the center $C$ of $M_n(R)$ is a subring of the ring $M_n(R)$ which is isomorphic to the ring $R$.
Thus ``$\mathcal{C}_1$ is isomorphic to $\mathcal{C}_2$" implies that ``$\mathcal{R}_1$ is isomorphic to $\mathcal{R}_2$".
We have then proved that $\mathcal{M}_1$ is isomorphic to $\mathcal{M}_2$ if and only if $\mathcal{R}_1$ is isomorphic to $\mathcal{R}_2$.
The result follows then from Theorem \ref{ind2}.
\ep

\hs We look now for similar results for groups. If $R$ is a ring with unity  element, i.e. a unitary ring,  then the set
$GL_n(R)$ of invertible matrices of $M_n(R)$ is a group for the multiplication
of matrices. It is first order interpretable in the ring $R$ because it is the set of matrices $M \in M_n(R)$
such that the determinant $det(M)$ of $M$ is invertible in $R$, see \cite{NiesBSL}.
This implies that if $R$ is an $\om$-tree-automatic unitary  ring   then the group $GL_n(R)$ is also $\om$-tree-automatic.

\hs Another interesting group is the unitriangular group $UT_n(R)$ for some integer $n\geq 3$ and $R$ a unitary ring.
A matrix $M\in M_n(R)$ is in the  group  $UT_n(R)$  if and only if it is an upper triangular matrix
which has only coefficients $1$ on the diagonal, where $1$ is the unity element for the second operation of $R$. The group $UT_n(R)$ is also first order
interpretable in the ring $R$ and it is a subgroup of the group $GL_n(R)$. We recall now the classical notion of nilpotent group.
The center of a group $G$ is denoted $Z(G)$. A group $G$ is said to be nilpotent of class $1$ iff $G$ is non-trivial and abelian. The group $G$
is said to be nilpotent of class $c+1$ if and only if the group $G / Z(G)$ is  nilpotent of class $c$.
The group $UT_n(R)$ is a classical example of nilpotent group of class $n-1$, see \cite{Belegradek}.

\hs We denote
$\mathcal{U}_{i,n}=UT_n(\mathcal{R}_i)$ for each $i \in \{1, 2\}$.
We have seen that the groups   $\mathcal{U}_{i,n}$ are first order interpretable in the ring $\mathcal{R}_i$. Thus the groups
 $\mathcal{U}_{i,n}$  are
$\om$-tree-automatic. We can now state the following result.

\begin{The}\label{ind4}
\noi
\begin{enumerate}
\ite  {\rm  (ZFC + CH)}  For each integer $n\geq 3$,
$\mathcal{U}_{1,n}$ is isomorphic to $\mathcal{U}_{2,n}$. \ite
{\rm   (ZFC + OCA) }   For each integer $n\geq 3$,
 $\mathcal{U}_{1,n}$ is not  isomorphic to $\mathcal{U}_{2,n}$ .
\end{enumerate}
\end{The}

\proo If $\mathcal{R}_1$ is isomorphic to $\mathcal{R}_2$ then it is clear that for each integer $n\geq 2$,
$\mathcal{U}_{1,n}$ is isomorphic to $\mathcal{U}_{2,n}$.
On the other hand, Belegradek  proved in \cite{Belegradek} that if  $UT_n(R)$ and $UT_n(S)$ are
isomorphic, for some integer $n\geq 3$ and some
commutative rings $R$ and $S$,  then the rings $R$ and $S$  are also  isomorphic.
Thus if for some integer $n\geq 3$ the groups $\mathcal{U}_{1,n}$ and  $\mathcal{U}_{2,n}$
 are isomorphic then $\mathcal{R}_1$ is isomorphic to $\mathcal{R}_2$.

\hs Then we have proved that, for each   integer $n\geq 3$,      $\mathcal{R}_1$ is isomorphic to $\mathcal{R}_2$  if and only if
the groups $\mathcal{U}_{1,n}$ and  $\mathcal{U}_{2,n}$ are isomorphic.
The  result follows then directly from Theorem \ref{ind2}.
\ep

\hs Then we can now infer  the following result.

\begin{Cor}
The isomorphism relation for $\om$-tree-automatic commutative rings (respectively,  non-commutative rings,
 groups,  nilpotent groups of class $n \geq 2$)
is not determined by the axiomatic system  {\rm  ZFC}.
\end{Cor}

\proo The result follows from Theorems \ref{ind2},  \ref{ind3},  \ref{ind4}, and the fact that the commutative rings $\mathcal{R}_i$, the non-commutative
rings $\mathcal{M}_i$, the groups $\mathcal{U}_{i,n}$, are all $\om$-tree-automatic structures.
\ep

\begin{Rem}{\rm
In  \cite{HjorthKMN08} the authors show that the isomorphism
relation for $\om$-automatic structures is not determined by the
axiomatic system {\rm  ZF}. They considered the two
$\om$-automatic groups $(\mathbb{R}, +)$ and $(\mathbb{R}, +)
\times (\mathbb{R}, +)$.  Assuming the axiom of choice is
satisfied,
 these two groups are isomorphic because they are both $\mathbb{Q}$-vectorial spaces of the same dimension
$2^{\aleph_0}$. But in Shelah's model of    {\rm  ZF}  where every
set of reals is Baire measurable the two groups are not
isomorphic, see \cite{HjorthKMN08}. 
We have then proved here a
stronger result in the case of $\om$-tree-automatic structures:
the isomorphism relation is not determined by the axiomatic system
{\rm  ZFC} (and not only {\rm  ZF}). Moreover we have proved our
result not only for the class of all $\om$-tree-automatic
structures and  for the class of $\om$-tree-automatic groups, but
also for the classes of $\om$-tree-automatic boolean algebras
(respectively, commutative rings,  non-commutative rings,
 nilpotent groups of class $n \geq 2$). }
\end{Rem}

\begin{Rem}{\rm
The two boolean algebras $\mathcal{B}_1$ and  $\mathcal{B}_2$ are
not isomorphic in every model of  {\rm  ZFC}. On the other hand,
they are always elementarily equivalent
 because they are two atomless boolean algebras and the first order theory of atomless boolean algebras is complete. In a similar way the rings
$\mathcal{R}_1$ and  $\mathcal{R}_2$ are
always elementarily equivalent. This implies that the rings $\mathcal{M}_1$ and  $\mathcal{M}_2$ (respectively, the groups
$\mathcal{U}_{1,n}$ and $\mathcal{U}_{2,n}$,
for some integer $n\geq 3$) are also always elementarily equivalent because the ring   $M_n(R)$ (respectively, the group $U_n(R)$ for $n\geq 3$)
is  first order interpretable in the ring $R$
without parameters and {\bf uniformly in the ring} $R$, see \cite[page 20]{Belegradek}. }

\end{Rem}

\noi An $\om$-tree-automatic presentation of a structure is given by a tuple of Muller tree automata
$(\mathcal{A}, \mathcal{A}_=,  (\mathcal{A}_i)_{1 \leq i \leq n})$, where  $L(\mathcal{A})\subseteq T_\Si^\om$, the automaton $\mathcal{A}_=$ accepts
an equivalence relation $E_\equiv $  on $L(\mathcal{A})$,  and for each $i \in [1, n]$, the automaton $\mathcal{A}_i$ accepts an $n_i$-ary relation $R_i$ on
$L(\mathcal{A})$ such that $E_\equiv$ is compatible with $R_i$. Then the tuple of automata
$(\mathcal{A}, \mathcal{A}_=,  (\mathcal{A}_i)_{1 \leq i \leq n})$ gives an $\om$-tree-automatic presentation
of the quotient  structure $( L(\mathcal{A}),  (R_i)_{1 \leq i \leq n}) ) / E_\equiv$.
\nl Notice that the tuple of automata can be coded by a finite sequence of symbols, hence by a unique integer $N$.
If $N$ is the code of the tuple of  Muller tree automata $(\mathcal{A}, \mathcal{A}_=,  (\mathcal{A}_i)_{1 \leq i \leq n})$
 we shall denote $\mathcal{S}_N$ the
$\om$-tree-automatic structure $( L(\mathcal{A}),  (R_i)_{1 \leq i \leq n}) ) / E_\equiv$.

\hs
The isomorphism  problem  for $\om$-tree-automatic structures is:

$$\{ (n, m) \in \mathbb{N}^2  \mid \mathcal{S}_n \mbox{ is isomorphic to } \mathcal{S}_m \}.$$

\noi A similar definition is given for automatic structures presentable by finite automata reading finite words, or
for $\om$-automatic structures presentable by B\"uchi automata reading infinite words.
It is proved in \cite{KNRS} that the isomorphism problem for automatic structures is $\Si_1^1$-complete.
The authors proved in   \cite{HjorthKMN08}
that the  isomorphism problem for $\om$-automatic structures is not
a $\Si_2^1$-set. In fact their proof implies also  that this isomorphism problem is not a $\Pi_2^1$-set.
Moreover  this is also the case for the restricted class of  $\om$-automatic (abelian) groups and   for the class of all  $\om$-tree-automatic structures
which is an extension of the class of $\om$-automatic structures.

\hs We can now  infer from above independence results  some similar results for many other classes of
$\om$-tree-automatic structures.

\begin{The}
The  isomorphism problem for  $\om$-tree-automatic boolean algebras (respectively, partial orders, rings, commutative rings,
non commutative rings, non commutative groups, nilpotent  groups of class $n\geq 2$)
is neither a $\Si_2^1$-set nor a  $\Pi_2^1$-set.
\end{The}

\proo We prove first the result for $\om$-tree-automatic boolean
algebras.  By Theorem \ref{ind} we know that if  {\rm ZF} hence
also {\rm ZFC} is consistent then there is a model ${\bf  V}$  of
{\rm (ZFC + OCA) } in which the boolean algebra $\mathcal{B}_1$ is
not isomorphic to the boolean algebra $\mathcal{B}_2$. But the
inner model ${\bf L }$  of constructible sets in  ${\bf  V}$  is a
model of {\rm  (ZFC + CH)} so in this model the two boolean
algebras $\mathcal{B}_1$ and $\mathcal{B}_2$ are  isomorphic. We
have also proved that these  two boolean algebras are
$\om$-tree-automatic.

\hs On the other hand, Schoenfield's Absoluteness Theorem implies that every $\Si_2^1$-set (respectively,  $\Pi_2^1$-set)
 is absolute for all inner models of {\rm  (ZF + DC)}, where {\rm  (DC)} is a weak version of the axiom of choice called axiom of dependent choice which
holds in particular in the inner model ${\bf L }$,
see \cite[page 490]{Jech}.

\hs In particular, if the isomorphism  problem  for  $\om$-tree-automatic boolean algebras was a $\Si_2^1$-set
(respectively,  a $\Pi_2^1$-set), then it could not be a different subset of $\mathbb{N}^2$ in the models  ${\bf  V}$  and   ${\bf L }$ considered above.
Thus the isomorphism  problem  for  $\om$-tree-automatic boolean algebras is neither a  $\Si_2^1$-set nor a $\Pi_2^1$-set.

\hs The other cases of  partial orders, rings, commutative rings,
non commutative rings, non commutative groups, nilpotent  groups of class $n\geq 2$, follow in the same way from
Theorems \ref{ind2}, \ref{ind3}, and \ref{ind4}.
\ep

\begin{Rem}{\rm
The set of codes of tuples of Muller tree automata which form
$\om$-tree-automatic presentations of boolean algebras is recursive because the first-order theory of  boolean algebras is finitely axiomatizable and
the first-order theory of an $\om$-tree-automatic structure is decidable by Theorem \ref{dec}. The same result holds in the cases of
partial orders (respectively, rings, commutative rings,
non commutative rings, non commutative groups).
Moreover this is also the case for  nilpotent  groups of class $n\geq 2$. Indeed it is decidable whether a tuple of Muller tree automata
$(\mathcal{A}, \mathcal{A}_=,  \mathcal{A}_1)$ is an $\om$-tree-automatic presentation of an abelian or a non-abelian group $(G, +)$. If it is an
$\om$-tree-automatic presentation of a non-abelian group $(G, +)$, then the center $Z(G)$ is first-order definable hence represented by a regular subset
of $L(\mathcal{A})$, and we can get an $\om$-tree-automatic presentation of the quotient  group $G / Z(G)$. If the group $G / Z(G)$ is
non-trivial and abelian then the group $G$ is   nilpotent of class $2$, and this can be determined again from Theorem \ref{dec}.  If  $G / Z(G)$ is
not abelian  then we can iterate the process, construct an $\om$-tree-automatic presentation of the quotient  group
$( G / Z(G) )  / Z(G / Z(G))$ and decide whether
this group is (non-trivial and) abelian. If this is the case then $G$ is nilpotent of class $3$, and so on.
}
\end{Rem}

\hs {\bf Acknowledgements.}
We are indebted to Anatole Khélif  who was generously sharing with  us his expertise about unitriangular groups 
and to the anonymous referees for very useful comments 
on a preliminary version of this paper.

\end{document}